\documentclass{amsart}
\usepackage{amssymb}
\usepackage{ifthen}
\usepackage{graphicx}

\newtheorem{thm}{Theorem}
\newtheorem{cor}{Corollary}
\newtheorem{lem}{Lemma}

\newtheorem{rem}{Remark}

\newtheorem{conj}{Conjecture}
\theoremstyle{definition}
\newtheorem{example}[equation]{Example}
\newtheorem{prob}[equation]{Problem}

\newcounter {own}
\def\theown {\thesection       .\arabic{own}}

\newenvironment{pf}[1][]{%
 \vskip 3mm
 \noindent
 \ifthenelse{\equal{#1}{}}%
  {{\slshape Proof. }}%
  {{\slshape #1.} }%
 }%
{\qed\bigskip}

\newcounter{alphabet}
\newcounter{tmp}

\newcommand{\ID}{{\mathbb D}}

\newcommand{\D}{{\mathbb D}}

\newcommand{\real}{{\operatorname{Re}\,}}

\def\be{\begin{equation}}
\def\ee{\end{equation}}

\newcommand{\bee}{\begin{enumerate}}
\newcommand{\eee}{\end{enumerate}}

\newcommand{\blem}{\begin{lem}}
\newcommand{\elem}{\end{lem}}
\newcommand{\bthm}{\begin{thm}}
\newcommand{\ethm}{\end{thm}}
\newcommand{\bcor}{\begin{cor}}
\newcommand{\ecor}{\end{cor}}
\newcommand{\beg}{\begin{example}}
\newcommand{\eeg}{\end{example}}
\newcommand{\begs}{\begin{examples}}
\newcommand{\eegs}{\end{examples}}
\newcommand{\bdefe}{\begin{defin}}
\newcommand{\edefe}{\end{defin}}
\newcommand{\bprob}{\begin{prob}}
\newcommand{\eprob}{\end{prob}}
\newcommand{\bei}{\begin{itemize}}
\newcommand{\eei}{\end{itemize}}

\newcommand{\bcon}{\begin{conj}}
\newcommand{\econ}{\end{conj}}
\newcommand{\bcons}{\begin{conjs}}
\newcommand{\econs}{\end{conjs}}
\newcommand{\bprop}{\begin{propo}}
\newcommand{\eprop}{\end{propo}}
\newcommand{\br}{\begin{rem}}
\newcommand{\er}{\end{rem}}
\newcommand{\brs}{\begin{rems}}
\newcommand{\ers}{\end{rems}}
\newcommand{\bo}{\begin{obser}}
\newcommand{\eo}{\end{obser}}
\newcommand{\bos}{\begin{obsers}}
\newcommand{\eos}{\end{obsers}}
\newcommand{\bpf}{\begin{pf}}
\newcommand{\epf}{\end{pf}}
\newcommand{\ba}{\begin{array}}
\newcommand{\ea}{\end{array}}
\newcommand{\beq}{\begin{eqnarray}}
\newcommand{\beqq}{\begin{eqnarray*}}
\newcommand{\eeq}{\end{eqnarray}}
\newcommand{\eeqq}{\end{eqnarray*}}

\begin{document}
\bibliographystyle{amsplain}

\title[A class of univalent functions with real coefficients]{A class of univalent functions with real coefficients}

\author[M. Obradovi\'{c}]{Milutin Obradovi\'{c}}
\address{Department of Mathematics,
Faculty of Civil Engineering, University of Belgrade,
Bulevar Kralja Aleksandra 73, 11000, Belgrade, Serbia}
\email{obrad@grf.bg.ac.rs}

\author[N. Tuneski]{Nikola Tuneski}
\address{Department of Mathematics and Informatics, Faculty of Mechanical Engineering, Ss. Cyril and Methodius
University in Skopje, Karpo\v{s} II b.b., 1000 Skopje, Republic of Macedonia.}
\email{nikola.tuneski@mf.edu.mk}







\maketitle
\pagestyle{myheadings}

\begin{abstract}
In this paper we study class $\mathcal{S}^+$ of univalent functions $f$ such that $\frac{z}{f(z)}$ has real and positive coefficients. For such functions we give estimates of the Fekete-Szeg\H{o} functional and sharp estimates of their initial coefficients and logarithmic coefficients. Also, we present necessary and sufficient conditions for $f\in \mathcal{S}^+$ to be starlike of order $1/2$.
\end{abstract}


\section{Introduction}
Let $\mathcal{A}$ be the class of functions $f$ that are analytic  in the open unit disc $\D=\{z:|z|<1\}$ of the form $f(z)=z+a_2z^2+a_3z^3+\cdots$. Then the class of starlike functions of order $\alpha$, $0\le\alpha<1$, is defined by
\[ \mathcal{S}^*(\alpha)=\left\{ f\in\mathcal{A} : \real \frac{zf'(z)}{f(z)}>\alpha, \, z\in\D\right\}, \]
while $\mathcal{S}^*\equiv \mathcal{S}^*(0)$ is the well known class of starlike functions mapping the unit disc onto a starlike region $D$, i.e.,
\[w \in f(D) \quad \Leftrightarrow \quad tw\in f(\D) \mbox{ for all } t\in[0,1].  \]
More on this classes can be found in \cite{TTV} and \cite{duren}.

\medskip

Further, let $\mathcal{S}^{+}$ denote the class of univalent functions in the unit disc
with the next representation
\be\label{eq 99}
\frac{z}{f(z)}=1+b_{1}z+b_{2}z^{2}+\ldots, \quad b_{n}\geq 0,\,n=1,2,3,\ldots
\ee
For example, the Silverman class (the class with negative coefficients) is included
in the class $\mathcal{S}^{+}$. Namely, that class consists of univalent functions of the form
$$f(z)=z-a_{2}z^{2}-a_{3}z^{3}-\ldots,\,\, a_{n}\geq 0,\,\, n=2,3,\ldots$$
which implies that
$$ \frac{z}{f(z)}= \frac{1}{1-a_{2}z-a_{3}z^{2}-\cdots},$$
i.e., $ \frac{z}{f(z)}$  has the form \eqref{eq 99}. Also, the Koebe function $k(z)=\frac{z}{(1+z)^{2}}\in\mathcal{S}^{+}$.
The next characterization is valid for the class $\mathcal{S}^{+}$ (see \cite{OP_2009}):
\be\label{eq 100}
f\in \mathcal{S}^{+} \quad\quad\Leftrightarrow \quad\quad\sum _{n=2}^{\infty}(n-1)b_{n}\leq 1.
\ee
From the relations \eqref{eq 99} and \eqref{eq 100} we have that
\be\label{eq 101}
b_{2}+2b_{3}\leq 1 \quad\quad (\Rightarrow \quad 0\leq b_{2}\leq 1,\quad 0\leq b_{3}\leq 1/2).
\ee
If we put  $f(z)=z+a_{2}z^{2}+\ldots$, then by using \eqref{eq 99} we easily obtain that
\be\label{eq 102}
 b_{1}=-a_{2},\quad b_{2}=a_{2}^{2}-a_{3}.
\ee
This implies that  $0 \leq b_{1}\leq 2$.
From \eqref{eq 99} we obtain
$$\log \frac{f(z)}{z}=-\log(1+b_{1}z+b_{2}z^{2}+\cdots) ,$$
or
$$\sum _{n=1}^{\infty}2\gamma_{n}z^{n}=-b_{1}z + \left(\frac{1}{2}b_{1}^{2}-b_{2}\right)z^{2}
+\left(-\frac{1}{3}b_{1}^{3}+b_{1}b_{2}-b_{3}\right)z^{3}+\cdots  .$$
(We call $\gamma_{n}, n=1,2,...$ the logarithmic coefficients of the function $f$.)
From the last relation we have
\be\label{eq 103}
 \left\{\begin{array}{l}
2\gamma_{1}=-b_{1}, \\
2\gamma _{2}=\frac{1}{2}b_{1}^{2}-b_{2} , \\
2\gamma _{3}= -\frac{1}{3}b_{1}^{3}+b_{1}b_{2}-b_{3}.
\end{array}
\right.
\ee

\medskip

For functions $f$ in $\mathcal{S}^+$ we give sharp estimates of their logarithmic coefficients $\gamma_1$, $\gamma_2$ and $\gamma_3$ of $f$ and lower and upper bound of the Fekete-Szeg\H{o} functional ($a_{3}-\gamma a_{2}^{2}$). Additionally, sharp estimates of coefficients $a_2$, $a_3$, $a_4$ and $a_5$ for functions in a class containing  $\mathcal{S}^+$ is given. At the end the relation between the class $\mathcal{S}^+$ and the class of starlike functions is studied.

\medskip

\section{Results over the coefficients}

\medskip

We start the section with a study of the Fekete-Szeg\H{o} functional for the functions in the class $\mathcal{S}^+$.

\bthm\label{13-th 27} For each $f\in \mathcal{S}^{+}$ we have
$$-1\leq a_{3}-\gamma a_{2}^{2}\leq \left\{\begin{array}{rl}
 1+2e^{-2\gamma/(1-\gamma)}, & 0\leq \gamma \leq \frac{\nu_{0}}{1+\nu_{0}}=0.456278\ldots \\
2(1-\gamma)\frac{(\nu_{0}+1)^{2}}{2\nu_{0}+1}, & \frac{\nu_{0}}{1+\nu_{0}}\leq\gamma < 1.
\end{array}
\right.$$
where $\nu_{0}=0.83927\ldots$  is the positive real root of the equation
\be\label{eq 104}
2(2\nu+1)e^{-2\nu}=1.
\ee
The lower bound is sharp due to the function $f_1(z)=\frac{z}{1+z^2}$.
\ethm

\medskip

\begin{proof} We will use the same method as in the proof of Fekete-Szego theorem for the class $\mathcal{S}$
(see \cite[Theorem 3.8, p. 104]{duren}).  First, from the relation \eqref{eq 102} we have that
$$ -1\leq a_{3}-a_{2}^{2}=-b_{2}\leq 0.$$
Since $a_{2}$ and $ a_{3}$ are real we can put (as in that proof)
$a_{2}=-2\int_{0}^{\infty}\varphi(t)dt$, where $\varphi$ is real function and $|\varphi(t)|\leq e^{-t}.$
If we put
$$\int_{0}^{\infty}[\varphi(t)]^{2}dt=\left(\nu+\frac{1}{2}\right)e^{-2\nu}, \quad 0\le\lambda<\infty,$$
then by Valiron-Landau lemma we have that $|a_{2}|\leq 2(\nu +1)e^{-\nu}$. Also, we have (as \cite{duren}):
$$ a_{3}- a_{2}^{2}=1-4\int_{0}^{\infty}[\varphi(t)]^{2}dt=1-4\left(\nu+\frac{1}{2}\right)e^{-2\nu}\leq 0$$
if and only if $0\leq \nu\leq\nu_{0}$, where $\nu_{0}=0.83927\ldots$ is the root of the equation
 \eqref{eq 104}.

Now, for $0\leq \gamma < 1$ and for $0\leq \nu\leq\nu_{0}$ we have that
 \beqq
 \begin{split}
a_{3}-\gamma a_{2}^{2}&\leq 4(1-\gamma)\left( \int_{0}^{\infty}\varphi(t)dt \right)^{2}
 -4 \int_{0}^{\infty}[\varphi(t)]^{2}dt +1 \\
 &= 4e^{-2\nu}\left[(1-\gamma)(\nu +1)^{2}-\left(\nu+\frac{1}{2}\right)\right]+1 \\
 &=: \psi (\nu).
 \end{split}
\eeqq
It is an elementary fact that the function $\psi$ has its maximum $\psi (\gamma/(1-\gamma))$ if
$\frac{\gamma}{1-\gamma}\in [0,\nu_{0}] $ and $\psi(\nu_{0})$ if $\frac{\gamma}{1-\gamma}\not\in [0,\nu_{0}]$,
which gives the right estimation in the theorem (in the second case we used that $\nu_{0}$ satisfies the equation
\eqref{eq 104}).

\medskip

On the other hand side
$ a_{3}-\gamma a_{2}^{2}\geq  a_{3}- a_{2}^{2}\geq -1 $.
\end{proof}

\medskip

Next we give sharp estimates of the first three logarithmic coefficients for functions in $\mathcal{S}^+$.

\bthm\label{13-th 28} Let $f\in \mathcal{S}^{+}$ and let $ \gamma _{1}, \gamma _{2}, \gamma _{3}$
be its logarithmic coefficients. Then
\begin{itemize}
  \item[$(a)$] $ -1\leq \gamma_{1}\leq 0$;
  \item[$(b)$] $ -\frac{1}{2}\leq\gamma_{2}\leq \frac{(\nu_{0}+1)^{2}}{2(2\nu_{0}+1)}=0.631464\ldots$, \\
where $\nu_{0}=0.83927\ldots$  is the solution of the equation
\eqref{eq 104};
\item[$(c)$] $-\frac{1}{4}\leq \gamma_{3}\leq \frac{1}{3}$.
\end{itemize}
Some of these results are the best possible.
\ethm

\begin{proof} $ $
\begin{itemize}
  \item[$(a)$] It is evident since $\gamma_{1}=-\frac{1}{2}b_{1}$ (from \eqref{eq 103}) and
$0 \leq b_{1}\leq 2$. The functions $f_{1}(z)=\frac{z}{1+z^{2}}$ and $f_{2}(z)=\frac{z}{(1+z)^{2}}$
show that the result is the best possible.
  \item[$(b)$] From \eqref{eq 102} and \eqref{eq 103} we have that
$$\gamma _{2}=\frac{1}{2}\left(\frac{1}{2}b_{1}^{2}-b_{2}\right)=\frac{1}{2}\left(a_{3}-\frac{1}{2}a_{2}^{2}\right)$$
and the result directly follows from Theorem \ref{13-th 27} for $\gamma=\frac{1}{2}$. For the function $f_{1}(z)=\frac{z}{1+z^{2}}$
we have that $\log\frac{f_{1}(z)}{z}=-\log(1+z^{2})=-z^{2}+\ldots$,
which means that left hand side  estimate is the best possible.

We were not able to prove sharpness of the right hand side of the inequality (the upper bound of $\gamma_2$), but it is worth pointing that the estimate goes in a line with the sharp estimate coresponding to the univalent functions, known to be (see \cite[Theorem 3.8]{duren} or \cite[p.136]{TTV})
\[ |\gamma_2| \le\frac12(1+2e^{-2}) = 0.635\ldots . \]

\item[$(c)$]
From \eqref{eq 103} we have
$$2\gamma _{3}= -\frac{1}{3}b_{1}^{3}+b_{1}b_{2}-b_{3}=:u(b_{1}),$$
where $$u(t)=-\frac{1}{3}t^{3}+b_{2}t-b_{3},\quad 0\leq t\leq 2.$$
Since $u'(t)=-t^{2}+b_{2}$ and $u'(t)=0$ for $t_{0}=\sqrt{b_{2}}$, then the function $u$
attains its maximum
$$u(t_{0})=u\left(\sqrt{b_{2}}\right)=\frac{2}{3}b_{2}^{3/2}-b_{3}
\leq \frac{2}{3}(1-2b_{3})^{3/2}-b_{3}\leq \frac{2}{3}, $$
because $b_{2}\leq 1-2b_{3}$ (see \eqref{eq 101}) and the last function is a decreasing
function of $b_{3},\,0\leq b_{3}\leq \frac{1}{2}$. This provides that $\gamma_{3}\leq \frac{1}{3}$.
For the function $f_{3}(z)=\frac{z}{1+z+z^{2}}$ we have
$$\log \frac{f_{3}(z)}{z}=-\log(1+z+z^{2})=-z-\frac{1}{2}z^{2}+\frac{2}{3}z^{3}+\cdots,$$
i.e. $\gamma_{3}=\frac{1}{3}$.

\medskip

As for lower bound for $\gamma_{3}$, by using \eqref{eq 103} and \eqref{eq 102}, we have
\beqq
\begin{split}
-2\gamma _{3} &= \frac{1}{3}b_{1}^{3}-b_{1}b_{2}+b_{3}\\
&= \frac{1}{3}b_{1}^{3}-b_{1}(b_{1}^{2}-a_{3})+b_{3}\\
&=-\frac{2}{3}b_{1}^{3}+a_{3}b_{1}+b_{3}\\
&= v(b_{1}),
\end{split}
\eeqq
where $$v(t)=-\frac{2}{3}t^{3}+a_{3}t+b_{3},\quad \quad (0\leq t \leq 2).$$
From here we have
$$v'(t)=-2t^{2}+a_{3}.$$
If $a_{3}\leq 0$, then $v'(t)\leq 0$, and if $a_{3}>0$ then we can write
$$v'(t)=-2(b_{1}^{2}-a_{3})-a_{3}=-2b_{2}-a_{3}$$ and also we have $v'(t)<0$, since $0\leq b_{2}\leq1$.
It means that the function $v$ is a decreasing function, which gives that
$$-2\gamma_{3}\leq v(0)=b_{3}\leq \frac{1}{2},$$  i.e $\gamma_{3}\geq -\frac{1}{4}$.
The function $f_{4}(z)=\frac{z}{1+z^{3}/2}$ shows that the result is the best possible.
\end{itemize}
\end{proof}

\medskip

Let $\mathcal{U}(\lambda)$, $0<\lambda \leq1$, denote the class of functions $f\in \mathcal{A}$
which satisfy the condition
\be\label{eq 7}
\left |\left (\frac{z}{f(z)} \right )^{2}f'(z)-1\right | < \lambda \quad\quad (z\in \ID).
\ee
We put $\mathcal{U}(1)\equiv\mathcal{U}.$
More about classes $\mathcal{U}$ and $\mathcal{U}(\lambda )$ we can find in \cite{OP_2011}, \cite{TTV}, \cite{OP_2016} and \cite{OP_2018}.

\medskip

Let $\mathcal{U}^{+}(\lambda )$, $0< \lambda \leq 1$, denote the class of functions $f$ satisfy the conditions
\eqref{eq 99} and \eqref{eq 7}. By using \eqref{eq 100} we can conclude that
$\mathcal{U}^{+}(1)\equiv \mathcal{S}^{+}$. For example, the function
\be\label{eq 8}
\begin{split}
f_{\lambda}(z) &= \frac{z}{1+(1+\lambda)z+\lambda z^{2}}\\
&=z-(1+\lambda)z^{2}+(1+\lambda+\lambda^{2})z^{3}
-(1+\lambda +\lambda^{2}+\lambda^{3})z^{4}+\cdots.
\end{split}
\ee
belongs to the class $\mathcal{U}^{+}(\lambda )$ and it is extremal in many cases.

\medskip

Also, if $f\in\mathcal{U}^{+}(\lambda )$ and has the form
\eqref{eq 99} then by \eqref{eq 100}:
$$\sum _{n=2}^{\infty}(n-1)b_{n}\leq \lambda, $$
which implies the appropriate inequalites :
\be\label{eq 9}
0\leq b_{2}\leq \lambda, \quad b_{2}+2b_{3}\leq \lambda, \quad b_{2}+2b_{3}+3b_{4}\leq \lambda,\ldots.
\ee

\medskip

For the coefficients of functions from the class $\mathcal{U}^{+}(\lambda )$, the next theorem is valid.

\medskip

\bthm\label{13-th 2}
If $f(z)=z+a_{2}z^{2}+a_{3}z^{3}+\cdots$ belongs to the class  $\mathcal{U}^{+}(\lambda )$, $0< \lambda \leq 1$,
then we have
\beqq
\begin{split}
-(1+\lambda)&\leq a_{2}\leq0,\\
-\lambda &\leq a_{3}\leq 1+\lambda+\lambda^{2}, \\
-(1+\lambda+\lambda^{2}+\lambda^{4})&\leq a_{4}\leq \frac{4\lambda}{3}\sqrt{\frac{2\lambda}{3}},\\
a_{5} & \ge\left\{
\begin{array}{cc}
-\lambda/3, & 0<\lambda\leq 4/27\\
-9\lambda^{2}/4, & 4/27 \leq \lambda \leq1
\end{array}
\right..
\end{split}
\eeqq
All these inequalities are sharp.
\ethm


\begin{proof}
For $f(z)=z+a_{2}z^{2}+a_{3}z^{3}+\cdots$ and $f\in\mathcal{U}(\lambda )$, $0< \lambda \leq 1$ it is shown in \cite{OP_2018} the next sharp inequalities:
$$ |a_{2}|\leq 1+\lambda,\quad |a_{3}|\leq 1+\lambda+\lambda^{2},\quad |a_{4}|\leq 1+\lambda+\lambda^{2}+\lambda^{3}.$$
In the same paper the authors conjectured that $|a_{n}|\leq \sum_{k=0}^{n-1}\lambda^{k}$. Since the function $f_{\lambda}$ defined by \eqref{eq 8}
belongs to the class $\mathcal{U}^{+}(\lambda )$, then the lower bounds for $a_{2}$ and $a_{4}$ and the upper bounds for $a_{3}$ are valid and sharp. We only need to prove the lower bounds for $a_{3}$ and $a_{5}$ and the upper bounds for $a_{2}$ and $a_{4}$.

\medskip

If $f(z)=z+a_{2}z^{2}+a_{3}z^{3}+\cdots$ and $f$ has the form \eqref{eq 99},  then by comparing the coefficients we easily conclude that
\be\label{eq10}
\left\{\begin{array}{l}
a_{2}=-b_{1},\\
a_{3}=-b_{2}+b_{1}^{2}, \\
a_{4}=-b_{3}+2b_{1}b_{2}-b_{1}^{3},\\
a_{5}=-b_{4}+b_{2}^{2}+2b_{1}b_{3}-3b_{1}^{2}b_{2}+b_{1}^{4}.
\end{array}
\right.
\ee

\medskip

From $a_{2}=-b_{1}$ and $b_{1}\geq 0$, we have $a_{2}\leq 0 $. Also, by using
\eqref{eq 9} and \eqref{eq10}, we obtain
$$-a_{3}=b_{2}-b_{1}^{2}\leq b_{2}\leq \lambda,$$
which implies $a_{3}\geq -\lambda$. The function $f_{6}(z)=\frac{z}{1+\lambda z^{2}}(=z-\lambda z^{3}+\cdots)$
shows that two previous results are the best possible.

\medskip

Further, from \eqref{eq10} we have
$$a_{4}=-b_{3}+2b_{1}b_{2}-b_{1}^{3}\leq 2b_{2}b_{1}-b_{1}^{3}=:w(b_{1}),$$
where $0\leq b_{1}\leq 1+\lambda $ (since $b_{1}=-a_{2}\leq 1+\lambda $). It is an elementary
fact to get that the function $w$ has its maximun $\frac{4b_{2}}{3}\sqrt{\frac{2b_{2}}{3}}$ for
$b_{1}=\sqrt{\frac{2b_{2}}{3}}.$ It means that
$$a_{4}\leq \frac{4b_{2}}{3}\sqrt{\frac{2b_{2}}{3}}\leq \frac{4\lambda}{3}\sqrt{\frac{2\lambda}{3}},$$
since $0\leq b_{2}\leq \lambda.$ The function
$$f_{7}(z)=\frac{z}{1+\sqrt{\frac{2\lambda}{3}}z+\lambda z^{2}}$$
shows that the result is the best possible.

\medskip

Finally, from \eqref{eq10} we also have
\beqq
\begin{split}
-a_{5}&= b_{4}-b_{2}^{2}-2b_{1}b_{3}+3b_{1}^{2}b_{2}-b_{1}^{4}\\
&\leq b_{4}+3b_{1}^{2}b_{2}-b_{1}^{4}\\
&\leq \frac{1}{3}(\lambda-b_{2})+3b_{2}b_{1}^{2}-b_{1}^{4}\\
&\leq \frac{9}{4}b_{2}^{2} +\frac13(\lambda-b_2)\\
&\leq
\left\{
\begin{array}{cc}
\lambda/3, & 0<\lambda\leq4/27\\
9\lambda^{2}/4 & 4/27 \leq \lambda \leq1
\end{array}
\right.,
\end{split}
\eeqq
where we used the relation \eqref{eq 9} and the same method as in the previous case. The functions
$$f_{2}(z)=\frac{z}{1+\sqrt{\frac{3\lambda}{2}}z+\lambda z^{2}} \quad \mbox{and} \quad f_8(z) = \frac{z}{1+\frac{\lambda}{3}z^4} $$
show that the result is the best possible.
\end{proof}

\medskip

For $\lambda=1$ from the previous theorem we have

\begin{cor}
Let $f(z)=z+a_{2}z^{2}+a_{3}z^{3}+\cdots$ belong to the class $\mathcal{S}^{+}$. Then we have the next sharp
inequalities
$$-2\leq a_{2}\leq0,\quad -1 \leq a_{3}\leq 3, \quad -4\leq a_{4}\leq \frac{4}{3}\sqrt{\frac{2}{3}}, \quad -\frac{9}{4}\leq  a_{5}\leq 5.$$
We note that upper bound for $a_{5}$ follows from de Brange's theorem.
\end{cor}

\medskip

\section{Relation with starlike functions}

\medskip

In this section we study the relation between the class $\mathcal{S}^+$ and the class of starlike functions.

\medskip

\bthm\label{13-th 29}
 Let $f\in\mathcal{A}$  and satisfy the condition \eqref{eq 99}.
 Then the condition
 \be\label{eq 105}
\sum_{n=1}^{\infty}(2n-1)b_{n}\leq 1
\ee
is  necessary and sufficient for $f$ to be in the class
$\mathcal{S}^{\star}(1/2)$.
\ethm

\medskip

\begin{proof}  The sufficient condition follows from the result given in the paper of Reade,
Silverman and Todorov \cite{Reade}.

\medskip

Let's prove the necessary case.
If $f\in {\mathcal S}^{\star}(\frac{1}{2})$, then
$$
\left|\frac{\frac{zf'(z)}{f(z)}-1}{\frac{zf'(z)}{f(z)}}\right|<1\quad\quad(z\in\D)
$$
 or equivalently
$$
\frac{\left|z\left(\frac{z}{f(z)}\right)'\right|}{\left|\frac{z}{f(z)}-z\left(\frac{z}{f(z)}\right)'\right|}<1 \quad\quad(z\in\D)
$$
and from here
$$
\frac{\left|\sum_{n=1}^{\infty}nb_{n}z^{n}\right|}{\left|1-\sum_{n=2}^{\infty}(n-1)b_{n}z^{n}\right|}<1 \quad\quad(z\in\D).
$$
If $z=r$ ($0<r<1$)  we have from the last inequality that
$$
\frac{\sum_{n=1}^{\infty}nb_nr^{n}}{1-\sum_{n=2}^{\infty}(n-1)b_nr^{n}}<1,
$$
which implies the condition
$$
\sum_{n=1}^{\infty}(2n-1)b_nr^{n}<1.
$$
Finally, when $r\rightarrow 1$ we have
$$
\sum_{n=1}^{\infty}(2n-1)b_n\leq1,
$$
i.e., the relation  \eqref{eq 105}.
\end{proof}

\medskip

\begin{rem}
Since the class of convex functions is the
subset of the class $S^{\star}(1/2)$, then if a function $f$
is convex and $$\frac{z}{f(z)}=1+b_1 z+b_{2}z^{2}+\dots $$ with
$b_n\geq0$ for  $n=1,2,\ldots$, we have
$$\sum_{n=1}^{\infty}(2n-1)b_n\leq 1.$$ The converse is not true.
Namely, for the function $$ f(z)=\frac{z}{1+\frac{1}{3}z^2},$$ we
have that $$ \frac{z}{f(z)}=1+\frac{1}{3}z^2$$ and
$$ \sum_{n=1}^{\infty}(2n-1)b_n=1,$$ but
%
%
$$
1+\frac{zf''(z)}{f'(z)}=
\frac{1-2z^2+\frac{1}{9}z^4}{1-\frac{1}{9}z^4}
< 0
$$
for $z=r(0<r<1) $ and $r$ close to 1.
\end{rem}

\medskip

\bthm\label{13-th 30}
Let $f\in \mathcal{S}^{+}$ and let $b_{1}=0$, then $f\in \mathcal{S}^{\star}$ .
\ethm


\begin{proof} Since $f\in \mathcal{S}^{+}$, then $\sum _{n=2}^{\infty}(n-1)b_{n}\leq 1$, and since $b_{1}=0$, then
also $\sum _{n=2}^{\infty}(n-1)b_{n}\leq 1=1-b_{1},$ which implies, by result of Reade et al. (\cite{Reade}) (see the previous sited result in Theorem \ref{13-th 29}),
that $f\in \mathcal{S}^{\star}$.

We note that if $b_{1}=0$, then $\real \frac{f(z)}{z}>\frac{1}{2}$ $(z\in\D)$ since
$$|z/f(z)-1|\leq |z|^{2} \sum_{n=2}^{\infty}b_{n} \leq \sum_{n=2}^{\infty}(n-1)b_{n}\leq |z|^{2}<1 \quad\quad (z\in\D).$$
But under the condition of this theorem we do not have that $f\in \mathcal{S}^{\star}(1/2)$. For example, for the function $f_{1}(z)=\frac{z}{1+z^{2}}$ we have $b_1=0$, but
$\sum _{n=1}^{\infty}(2n-1)b_{n}=3$, which means that $f_1\notin \mathcal{S}^{\star}(1/2)$ (by the previous theorem).
\end{proof}

\medskip

\bthm\label{13-th 31}
Let $f\in \mathcal{S}^{+}$.
Then the function
\be\label{eq 106}
g(z)=z+\frac{1}{2}\left(\frac{z}{f(z)}-1-b_{1}z \right)
\ee
is univalent in $\ID$. More precisely,  $\real g'(z)>0$ $(z\in\D)$, $g\in \mathcal{S}^{\star}$ and $g\in \mathcal{U}$.
\ethm

\medskip

\begin{proof}
It is well-known that if $f(z)=z+a_{2}z^{2}+a_{3}z^{3}+\cdots$ and $\sum_{n=2}^{\infty}n|a_{n}|\leq 1$, then $\real f'(z)>0$ $(z\in\D)$ and $f\in \mathcal{S}^{\star}$ with
$\left|\frac{zf'(z)}{f(z)}-1\right|<1$ $(z\in\D)$. It is easily to prove those statement (in the second case better to consider the form $|zf'(z)-f(z)|<|f(z)|$).

\medskip

By  \eqref{eq 106} we have
$$g(z)= z+\sum _{n=2}^{\infty}\frac{1}{2}b_{n}z^{n}.$$
Since $f\in \mathcal{S}^{+}$ implies $\sum _{n=2}^{\infty}(n-1)b_{n}\leq 1$
and since $\frac{n}{2(n-1)}\leq 1$  for $n \geq 2$, then
$$\sum _{n=2}^{\infty}n \left(\frac{1}{2}b_{n}\right)=\sum _{n=2}^{\infty}(n-1)b_{n}\frac{n}{2(n-1)}\leq \sum _{n=2}^{\infty}(n-1)b_{n}\leq 1.$$
By previous remarks we have  $\real g'(z)>0$ $(z\in\D)$ and $g\in \mathcal{S}^{\star}$.
Also, $g\in \mathcal{U}$ by the result given in  \cite{OP_2011}.
\end{proof}

\end{document}